\UseRawInputEncoding

\documentclass[12pt]{amsart}
\usepackage{latexsym, amssymb, amsmath, amscd}
 \usepackage{eucal}

    \newtheorem{rema}{Remark}[section]
    \newtheorem{propo}[rema]{Proposition}

   \newtheorem{theo}[rema]{Theorem}
   \newtheorem{def-theo}[rema]{Definition-Theorem}
 \newtheorem{conj}[rema]{Conjecture}
   \newtheorem{defi}[rema]{Definition}
    \newtheorem{lemma}[rema]{Lemma}
    \newtheorem{corol}[rema]{Corollary}
     \newtheorem{exam}[rema]{Example}
  
  \newtheorem{rmks}[rema]{Remarks}

    \newtheorem{prob}[rema]{Problem}

	\newcommand{\p}{\partial}

 \newcommand{\pf}{{\it Proof:}\hspace{2ex}}
 
 \newcommand{\epfv}{\hspace{1em}$\Box$\vspace{1em}}

\newcommand{\bC}{{\mathbb C}}
\newcommand{\bZ}{{\mathbb Z}}

\newcommand{\bQ}{{\mathbb Q}}

\newcommand{\bF}{{\mathbb F}}

\newcommand{\cA}{{\mathcal A}}

\newcommand{\rad}{{\frak r}}
\newcommand{\nil}{\mbox{\rm nil\,}}
\newcommand{\I}{{\operatorname I}}

\newcommand{\cB}{{\mathcal B}}
\newcommand{\cD}{{\mathcal D}}

\newcommand{\cE}{{\mathcal E}}

\newcommand{\kx}{K[x]}

\newcommand{\Ker}{\mtype{\rm Ker\,}}
\newcommand{\im}{\mtype{Im\,}}

\newcommand{\mtype}{\operatorname}
\newcommand{\Der}{{\cD er}}
\newcommand{\Eder}{{\cE der}}
\newcommand{\cEnd}{{\cE nd}}

\title[Some Open Problems on LF or LN ($\cE$-)Derivations]  
{Some Open Problems on Locally Finite or Locally Nilpotent Derivations 
and $\cE$-Derivations}
 
  \author{Wenhua Zhao}      
    \date{\today}
\address{Department of Mathematics, Illinois State University, Normal, IL 61761. Email: wzhao@ilstu.edu}

\begin{document}

\begin{abstract}
Let $R$ be a commutative ring and $\cA$ an $R$-algebra. 
An $R$-$\cE$-derivation of $\cA$ is an $R$-linear map of 
the form $\I-\phi$ for some $R$-algebra endomorphism $\phi$ of $\cA$, 
where $\I$ denotes the identity map of $\cA$. In this paper we discuss 
some open problems on whether or not the image of 
a locally finite $R$-derivation or $R$-$\cE$-derivation of 
$\cA$ is a Mathieu subspace \cite{GIC, MS} of $\cA$, and whether or not   
a locally nilpotent $R$-derivation or $R$-$\cE$-derivation of 
$\cA$ maps every ideal of $\cA$ to a Mathieu subspace of $\cA$. 
We propose and discuss two conjectures which state that 
both questions above have positive answers if the base ring $R$ is 
a field of characteristic zero. 
We give some examples to show the necessity of 
the conditions of the two conjectures, and discuss 
some positive cases known in the literature.  
We also show some cases of the two conjectures.  
In particular, both the conjectures are 
proved for locally finite or locally nilpotent algebraic derivations 
and $R$-$\cE$-derivations of integral domains of characteristic zero.     
\end{abstract}

\keywords{Mathieu subspaces (Mathieu-Zhao spaces); the LNED conjecture; the LFED conjecture; locally finite or locally nilpotent derivations and $\cE$-derivations; idempotents; the Singer-Wermer Theorem}
   
\subjclass[2000]{47B47, 08A35, 16W25, 16D99}

%47B47 	Commutators, derivations, elementary operators, etc.

% 16D99 None of the above, but in this section 	"Associative rings and algebras" {For the commutative case, see 13-XX}

% 13G99	None of the above, but in this section of "Integral domains".

%08A35 	Automorphisms, endomorphisms

%16Wxx Rings and algebras with additional structure

%16W20 Automorphisms and endomorphisms

%16W25 Derivations, actions of Lie algebras

%16N40 Nil and nilpotent radicals, sets, ideals, rings
%
%16D70 Structure and classification (except as in 16Gxx), direct sum decomposition, cancellation
%
%16D99 None of the above, but in this section

%16S34 Group rings [see also 20C05, 20C07], Laurent polynomial rings

%13N10: Rings of differential operators and their modules 
%14R15: Jacobian problem.
%33C45 Orthogonal polynomials and functions of hypergeometric type (Jacobi, Laguerre, Hermite, Askey scheme, etc.) [see also 42C05 for general orthogonal polynomials and functions]. 
%32C38 Sheaves of differential operators and their modules, $D$-modules. 
%32W99: Differential operators in several 
%    complex variables/None of the above, but in this section. 

\thanks{The author has been partially supported 
by the Simons Foundation grant 278638}

 \bibliographystyle{alpha}
    \maketitle

%\tableofcontents

\renewcommand{\theequation}{\thesection.\arabic{equation}}
\renewcommand{\therema}{\thesection.\arabic{rema}}
\setcounter{equation}{0}
\setcounter{rema}{0}
\setcounter{section}{0}

\section{\bf Introduction}\label{S1}
   
Let $R$ be a unital ring (not necessarily commutative) and $\cA$ an $R$-algebra. 
We denote by $1_\cA$ or simply $1$ the identity element of $\cA$, if $\cA$ is unital, and 
$\I_\cA$ or simply $\I$ the identity map of $\cA$, if $\cA$ is clear in the context.  

An $R$-linear endomorphism $\eta$ of $\cA$ is said to be {\it locally nilpotent} (LN) 
if for each $a\in \cA$ there exists $m\ge 1$ such that $\eta^m(a)=0$, 
and {\it locally finite} (LF) if for each $a\in \cA$ the $R$-submodule spanned 
by $\eta^i(a)$ $(i\ge 0)$ over $R$ is finitely generated.    

By an $R$-derivation $D$ of $\cA$ we mean an $R$-linear map 
$D:\cA \to \cA$ that satisfies $D(ab)=D(a)b+aD(b)$ for all $a, b\in \cA$. 
By an $R$-$\cE$-derivation $\delta$ of $\cA$ 
we mean an $R$-linear map $\delta:\cA \to \cA$ such that for all 
$a, b\in \cA$ the following equation holds:
\begin{align}\label{ProdRule2}
\delta(ab)=\delta(a)b+a\delta(b)-\delta(a)\delta(b).  
\end{align}

It is easy to verify that $\delta$ is an $R$-$\cE$-derivation of $\cA$, 
if and only if $\delta=\I-\phi$ for some $R$-algebra endomorphism 
$\phi$ of $\cA$. Therefore an $R$-$\cE$-derivation is  
a special so-called $(s_1, s_2)$-derivation introduced 
by N. Jacobson \cite{J} and also a special 
semi-derivation introduced 
by J. Bergen in \cite{Bergen}. $R$-$\cE$-derivations have also been 
studied by many others under some different names such as 
$f$-derivations in \cite{E0, E} and 
$\phi$-derivations in \cite{BFF, BV}, etc..

We denote by $\cEnd_R(\cA)$ the set of all 
$R$-algebra endomorphisms of $\cA$, $\Der_R(\cA)$ the set of all 
$R$-derivations of $\cA$, and $\Eder_R(\cA)$ the set of all 
$R$-$\cE$-derivations of $\cA$. Furthermore, for each $R$-linear endomorphism 
$\eta$ of $\cA$ we denote by $\im \eta$ the {\it image} of $\eta$, i.e., 
$\im\eta\!:=\eta(\cA)$, and $\Ker \eta$ the {\it kernel} of $\eta$.
When $\eta$ is an $R$-derivation or $R$-$\cE$-derivation, we also denote  
by $\cA^{{}^\eta}$ the kernel of $\eta$. 

For each $R$-derivation or $R$-$\cE$-derivation of an $R$-algebra $\cA$, 
it is easy to see that the kernel $\cA^\delta$ is an $R$-subalgebra. Actually, 
if $\delta=\I-\phi$ for some $\phi\in\cEnd_R(\cA)$, the kernel $\cA^\delta$ of $\delta$ 
coincides with the $R$-subalgebra of the elements of $\cA$ that are fixed by $\phi$.
The kernels of derivations as well as the kernels of $\cE$-derivations 
(i.e., the subalgebra fixed by algebra endomorphisms) are among the most 
studied subjects and play important roles in various different areas 
(e.g., see \cite{N}, \cite{F}, \cite{E} and the references therein).  

On the other hand, the images, especially, their possible algebraic structures, 
of derivations or $\cE$-derivations have barely been studied. It is presumably because that in general they  
are not even closed under the multiplication of the algebra. 
However, recent studies (e.g., see \cite{EWZ}, \cite{Idem}--\cite{OneVariableCase}) 
show that the images of certain derivations 
and $\cE$-derivations do possess some algebraic structure. 
To be more precise, we first need to recall the following 
notion introduced in \cite{GIC, MS}.    

\begin{defi}   \label{Def-MS}
Let $\vartheta$ represent the words:  
$\mtype{left}$, $\mtype{right}$, or $\mtype{two-sided}$. An $R$-subspace $V$ 
of an $R$-algebra $\cA$ is said to be a $\vartheta$-Mathieu subspace ($\vartheta$-MS)  
of $\cA$ if for all $a, b, c\in \cA$ with $a^m\in V$ for all $m\ge 1$,  
the following conditions hold: 
\begin{enumerate}
  \item[$1)$] $ba^m\in V$ for all $m\gg 0$, if $\vartheta=\mtype{left}$;
  \item[$2)$] $a^mc\in V$ for all $m\gg 0$, if $\vartheta=\mtype{right}$;
  %\item $ba^mc\in V$ for all $m\gg 0$, if $\vartheta=pre-two-sided$;
  \item[$3)$] $b a^m c \in V$ for all $m\gg 0$, if $\vartheta=\mtype{two-sided}$. 
\end{enumerate} 
\end{defi}

A two-sided MS will also be simply called a MS. For an arbitrary 
ring $\cB$, the $\vartheta$-MSs of $\cB$ are defined by viewing $\cB$ as 
an algebra over $\bZ$. Some more remarks on the notion of 
MS are as follows.

First, the introduction of the notion in \cite{GIC} and \cite{MS}   
is mainly motivated by the Mathieu conjecture in \cite{Ma} 
and the Image conjecture in \cite{IC}, both of which are 
motivated by and also imply the well-known Jacobian conjecture  
that was first proposed by O. H. Keller in $1939$ \cite{K}. 
See also \cite{BCW} and \cite{E}.  
But, a more interesting aspect of the new 
notion is that it provides a natural but highly non-trivial 
generalization of the corner-stone notion of ideals 
of associative algebras. 

Second, a Mathieu subspace is also called  a {\it Mathieu-Zhao space} 
in the literature (e.g., see \cite{DEZ, EN, EH}, etc.) 
as first suggested by A. van den Essen \cite{E2}.

Third, the following notion, first introduced in \cite{MS},  
is closely related with MSs, although it is defined for 
all $R$-subspaces, or even arbitrary subsets, 
of $R$-algebras. 

\begin{defi} \cite[p.\,247]{MS} \label{Def-rad}
Let $V$ be an $R$-subspace of an $R$-algebra $\cA$. 
We define the {\it radical} of $V$, 
denoted by $\rad(V)$, to be the set 
of $a\in \cA$ such that $a^m\in V$ for 
all $m\gg 0$.  
\end{defi}

When $\cA$ is commutative and $V$ is an ideal of $\cA$, $\rad(V)$ 
coincides with the radical of $V$. So this new notion is also 
interesting on its own right. 
It is also crucial for the study of MSs. For example, the following lemma 
can be easily verified, and will be frequently used (implicitly) in this paper.  

\begin{lemma}\label{I-MS-Lma}
Let $V$ be an $R$-subspace of an $R$-algebra $\cA$, and $I$ an ideal of $\cA$. 
If $I\subseteq V$ and $\rad(I)=\rad(V)$. Then $V$ is a MS of $\cA$.  
\end{lemma}

Now we propose the following problems on the image of derivations and $\cE$-derivations.

\begin{prob}[{\bf LFNED Problem}]\label{LFNED-Prob}
Let $R$ be a commutative base ring, $\cA$ an $R$-algebra and $\delta$ an $R$-derivation or $R$-$\cE$-derivation of $\cA$.  
\begin{enumerate}
\item[${\bf A)}$] Find the radical of $\delta(I)$ for all one-sided or 
two sided ideals $I$ of $\cA$.  
\item[${\bf B)}$] Decide which $R$-derivations and 
 $R$-$\cE$-derivations of $\cA$ have the image being a $\vartheta$-MS of $\cA$.  
   \item[${\bf C)}$] Decide which $R$-derivations and  
 $R$-$\cE$-derivations of $\cA$ map every $\vartheta$-ideal of $\cA$  
to a $\vartheta$-MS of $\cA$.   
\end{enumerate}
\end{prob}

Although the sufficient and necessary conditions for Problem $B)$ and $C)$ are currently far 
from being clear, based on the studies in \cite{EWZ}, \cite{Idem}--\cite{OneVariableCase} 
as well as some results that will be shown later in this paper, 
the following two conjectures 
seem to be more plausible.

\begin{conj}[{\bf The LFED Conjecture}]\label{LFED-Conj}
Let $K$ be a field of characteristic zero and $\cA$ a $K$-algebra. Then for every 
locally finite $K$-derivation or locally finite $K$-$\cE$-derivation $\delta$ of $\cA$, the image 
$\im \delta\!:=\delta(\cA)$ of $\delta$ is a (two-sided) MS of $\cA$.   
\end{conj}

\begin{conj}[{\bf The LNED Conjecture}]\label{LNED-Conj}
Let $K$ be a field of characteristic zero, $\cA$ a $K$-algebra 
and $\delta$ a locally nilpotent  $K$-derivation or a locally nilpotent $K$-$\cE$-derivation of $\cA$. 
Then for every $\vartheta$-ideal $I$ of $\cA$, the image 
$\delta(I)$ of $I$ under $\delta$ is a $\vartheta$-MS of $\cA$.   
\end{conj}
 
Throughout this paper we refer to the two conjectures above as 
{\it the LFED conjecture} and {\it the LNED conjecture}, respectively.

One motivation of the two conjectures above 
is that they may provide some new understandings 
on the LF or LN derivations and $\cE$-derivations. 
Another motivation is that they may produce many non-trivial 
examples of MSs, which will be beneficial and essential toward 
the further development of the desired theory of MSs.   

Two more remarks on the two conjectures above are as follows. 
Below we let $K$ be a field of characteristic zero 
and $\cA$ a $K$-algebra, unless stated otherwise. 

First, \label{LND=LNED} by van den Essen's one-to-one correspondence 
(see \cite{E0} or \cite[Proposition 2.1.3]{E}) between 
the set of LN $K$-derivations of $\cA$ and the set of LN $K$-$\cE$-derivations 
of $\cA$ and also \cite[Corollary 2.4]{Idem}, 
the LN $K$-derivation case and the LN $K$-$\cE$-derivation case 
of Conjecture \ref{LFED-Conj} are equivalent to each other.   
In other words,  Conjecture \ref{LFED-Conj} holds for all 
LN $K$-derivations of $\cA$, if and only if it holds for all 
LN $K$-$\cE$-derivations of $\cA$.
 
Second, for every $\vartheta$-MS\, $V$ of $\cA$ and idempotent $e\in V$ 
(i.e., $e^2=e$), by Definition \ref{Def-MS} it is easy to see that the principal 
$\vartheta$-ideal $(e)_\vartheta$ of $\cA$ generated by $e$ is contained in $V$.  
Therefore, we have the following weaker versions of 
Conjectures \ref{LFED-Conj} and \ref{LNED-Conj}.

\begin{conj}[{\bf The Idempotent Conjecture}]\label{Idem-Conj}
Let $K$ be a field of characteristic zero, $\cA$ a $K$-algebra 
and $\delta$ a $K$-derivation or $K$-$\cE$-derivation of $\cA$.
Then the following two statements hold: 
\begin{enumerate}
\item[${\bf A)}$] If $\delta$ is LF, then for all idempotents $e\in \im\delta$, the principle ideal $(e)$ is contained in $\im \delta$; 
\item[${\bf B)}$] If $\delta$ is LN, then for all $\vartheta$-ideal $I$ of $\cA$ and all idempotents $e\in \delta(I)$, the $\vartheta$-ideal $(e)_\vartheta$ is 
 contained in $\delta(I)$.    
\end{enumerate}
\end{conj}

Actually, if $\cA$ is algebraic over $K$, the statements 
$A)$ and $B)$ in Conjecture \ref{Idem-Conj}  are 
respectively equivalent to the LFED conjecture  and the LNED conjecture, 
due to the following characterization of MSs 
of $\cA$, which is a special case of \cite[Theorem $4.2$]{MS}.

\begin{theo} \label{IdemCri} 
Let $K$ be a field (of arbitrary characteristic)  
and $\cA$ a $K$-algebra that is algebraic over $K$. 
Then a $K$-subspace $V$ of $\cA$ is a $\vartheta$-MS of $\cA$, 
if and only if for every idempotent $e\in V$, the principal 
$\vartheta$-ideal $(e)_\vartheta$ of $\cA$ generated by $e$ is contained in $V$.
\end{theo} 

\vspace{3mm}

{\bf Arrangement}:  
In Section \ref{S2}, we mainly give some examples to show 
the necessity of the conditions in the LFED conjectures \ref{LFED-Conj} and 
the LNED conjectures \ref{LNED-Conj}. We also give some positive examples 
with certain weaker conditions. 
In Section \ref{S3}, we discuss some positive cases of 
Conjectures \ref{LFED-Conj},  \ref{LNED-Conj} and \ref{Idem-Conj}, 
which are either already known in the literature or 
can be derived from some other results in the literature. 

In Section \ref{S4}, we discuss the LFED conjectures \ref{LFED-Conj} 
in terms of the decompositions of the $K$-algebra $\cA$ associated with 
the Jordan-Chevalley decomposition of the LF $K$-derivations 
and $K$-$\cE$-derivations of $\cA$. Two other conjectures (see Conjectures \ref{Grading-Conj} and \ref{GeneralDK-Conj}) that are closely related with the semi-simple case of the LFED conjecture  
are also proposed and discussed.

In Section \ref{S5}, we show that the LFED conjecture \ref{LFED-Conj} holds 
for $\cE$-derivations associated with some special algebra endomorphisms such as 
projections and involutions, etc.. In Section  \ref{S6}, we study the LFED 
conjectures \ref{LFED-Conj} and the LNED conjectures \ref{LNED-Conj} for 
algebraic derivations and $\cE$-derivations of domains of characteristic zero. 
In particular, for integral domains $\cA$ of characteristic zero  
we show that both conjectures hold for LF or LN algebraic derivations and 
$\cE$-derivations of $\cA$ 
(see Proposition \ref{HDIntDomain} and Theorem \ref{AlgIntDomainThm}).    \\

{\bf Acknowledgment:} The author is very grateful to Professor Arno van den Essen 
for reading carefully an earlier version of the paper and pointing out some 
mistakes and typos, and in particular, 
for sending the author some counter-examples for an earlier  
(and stronger) version of Conjecture \ref{Grading-Conj}.

\renewcommand{\theequation}{\thesection.\arabic{equation}}
\renewcommand{\therema}{\thesection.\arabic{rema}}
\setcounter{equation}{0}
\setcounter{rema}{0}

\section{\bf Some Examples and Necessity of the Conditions of the LFED and LNED Conjectures}\label{S2}

In this section we give some examples to show that the conditions in 
the LFED conjecture \ref{LFED-Conj} and the LNED conjecture \ref{LNED-Conj} are necessary. 
We also give some (positive) examples with some weaker conditions. 

Throughout this section {\it $K$ stands for a field of characteristic zero.}
All the notations introduced in the previous section will also be in force.

First, the following two examples show that the LF (locally finite) condition is  
necessary for both the LFED and LNED conjectures. 

\begin{exam}\cite[Example $2.4$]{EWZ} Let $x$ and $y$ be two commutative free variables 
and $D = \p/\p x - y^2\p/\p y$. Then $D$ is not LF and 
$\im D$ is not a MS of the polynomial algebra $K[x, y]$.
\end{exam}

\begin{exam}\label{HD-CounterEx}
Let $x, y$ be two commutative variables and $\phi$ the $K$-algebra endomorphism of $K[x, y]$ 
such that $\phi(x)=x+1$ and $\phi(y)=y^2$. Set $\delta\!:=\I-\phi$. 
Then it is easy to see that $\delta$ is not LF, and 
$\im \delta \ne \bC[x, y]$, since each $f\in \im\delta$ 
with $\deg_y f\ge 1$ must have even degree in $y$. 
On the other hand, $1=\delta(-x)\in\im \delta$. 
Then it is easy to check (or by \cite[Lemma 4.5]{GIC}) that  
$\im\delta$ is not a MS of $K[x, y]$. 
\end{exam}

Next, the following two examples show that the LN (locally nilpotent) condition 
in the LNED conjecture \ref{LNED-Conj} is necessary and 
can not be replaced by the LF (locally finite) condition.

\begin{exam}\cite[Example $2.4$]{OneVariableCase} \label{xddxExample} %\label{xddxExample} there too.
Let $x$ be a free variable, $D=x\frac d{dx}$ and 
$I=(x^2-1)\kx$. Then $D$ is LF but the image $DI$ 
of $I$ under $D$ is not a MS of $K[x]$.
\end{exam}

\begin{exam}\cite[Example $3.6$]{OneVariableCase} \label{qxExample} %\label{qxExample} there too.
Let $K$, $x$, $I$ be as in Example \ref{xddxExample}, and 
$0\ne q\in K$ that is not a root unity. Let $\phi\in \cEnd_K(K[x])$ 
that maps $x$ to $qx$ and $\delta\!:=\I-\phi$. Then $\delta$ is LF but  
the image $\delta I$ of $I$ under $\delta$ is  
not a MS of $\kx$. 
\end{exam}

The following two examples show that the base field $K$ in 
the LFED and LNED Conjectures can not be replaced 
by a field of characteristic $p>0$.

\begin{exam}\cite[Example $2.7$]{IC}\label{p-CounterEx}
Let $\bF$ be a field of characteristic $p>0$, $x$ a free variable and 
$D\!:=d/ dx$. Then $D$ is LN but $\im D$ is not a MS of $\bF[x]$.  
\end{exam}

\begin{exam}\label{p-Count-LFED}
Let $\bF_2=\bZ/2\bZ$, $x$ a free variable, $\cA=\bF_2[x^{-1}, x]$, and 
$\phi\in \cEnd_{\bF_2}(\cA)$ that maps $x$ to $x^{-1}$. 
Then $\delta\!:=\I-\phi$ is LF but 
$\im\delta$ is not MS of $\cA$.
\end{exam}

\pf It is easy to see that $\delta$ is LF. 
To show the second statement,
let $u_m\!:=\delta(x^m)=x^m+x^{-m}$ for all $m\ge 1$, 
and $V$ be the $\bF_2$-subspace of $\cA$ spanned 
by $u_m$ $(m\ge 1)$. Then $V\subseteq \im \delta$. 
By the binomial formula (over $\bF_2$) and 
the fact $\binom{m}i=\binom{m}{m-i}$ for all $0\le i\le m$, 
it is easy to see that 
for all $k\ge 1$, we have $(x+x^{-1})^{2k-1} \in V$ and   
\begin{align*}
(x+x^{-1})^{2k} \equiv \binom{2k}k \mod V.
\end{align*}
Note that $\binom{2k}k$ is even, which can actually be seen by 
letting $x=1$ in the equation above. Therefore,  
$(x+x^{-1})^m\in V\subseteq \im\delta$ for all $m\ge 1$.

On the other hand, since $\phi^2=\I$, we have 
$\phi\delta=\delta$. Therefore, every $f\in \im\delta$ is fixed 
by $\phi$, i.e., $\phi(f)=f$. By this fact we see that    
$x(x+x^{-1})^m \not \in \im \delta$ for all $m\ge 1$.
Hence $\im\delta$ is not a MS of $\cA$.
\epfv

Although the LFED and LNED Conjectures can not be extended to all the algebras over  
a field of characteristic $p>0$ (as shown by the two examples above), 
the following example shows that the LFNED problem \ref{LFNED-Prob} 
is still interesting for some of these algebras. 

\begin{exam}\label{FrobeniCase}
Let $p\ge 2$ be a prime, $\bF_p=\bZ/p\bZ$, $x$ a free variable and $\phi$ 
the Frobenius endomorphism of $\bF_p[x]$, i.e., $\phi(f)=f^p$ for all $f\in \bF_p[x]$. 
Set $\delta\!:=\I-\phi$. Then $\rad\big(\im \delta\big)=\{0\}$. 
Consequently, $\delta$ maps every $\bF_p$-subspace of $\bF_p[x]$ to a MS 
of $\bF_p[x]$. 
\end{exam} 

Since $a^p=a$ for all $a\in \bF_p$, $\phi$ is actually an $\bF_p$-algebra endomorphism of 
$\bF_p[x]$. Then the conclusion of the example follows from 
\cite[Proposition 3.7]{OneVariableCase}. But, for the sake of completeness 
we include here a more straightforward proof. \\

\underline{\it Proof of Example \ref{FrobeniCase}}:\, We first show $1 \not \in \im \delta$. 
Otherwise, let $f\in \bF_p[x]$ such that
$\delta f=f-\phi(f)=f(x)-f^p(x)=1$.  
Then $\deg f =0$, i.e., $f \in \bF_p$. 
But in this case $f^p=f$, 
whence $f-\phi(f)=0$. Contradiction.  

Now assume $\rad(\im \delta)\ne \{0\}$ and let $0\ne f(x)\in \rad(\im \delta)$. 
Then $\deg f\ge 1$. Replacing $f$ by a power of $f$ we assume 
$f^m\in \im\delta$ for all $m\ge 1$.  
Let $h_m\in \bF_p [x]$ $(m\ge 1)$ such that    
\begin{align}\label{FrobeniCase-peq1}
f^m(x)=\delta h_m=h_m(x)-h_m^p(x).
\end{align}

If $f'(x)=\frac{df}{dx} (x)=0$, then 
$f(x)=\tilde f(x^p)$ for some 
$\tilde f(x)\in \bF_p[x]$. By the equation above with $m=1$ we have 
$h_1'(x)=0$, whence $h_1(x)=\tilde h_1(x^p)$ for some 
$\tilde h_1(x)\in \bF_p [x]$. The equation above with $m=1$ becomes
$\tilde f(x^p)=\tilde h_1(x^p)-\tilde h_1^p(x^p)$. Replacing 
$x^p$ by $x$ we have $\tilde f(x)=\tilde h_1(x)-\tilde h_1^p(x) 
\in \im\delta$. 

Applying the same arguments to $f^m$ $(m\ge 2)$ we see that there exists 
$\tilde f_m(x)\in \im\delta$ such that $\tilde f_m(x^p)=f^m(x)=\tilde f^m(x^p)$. 
Hence $ \tilde f_m(x)=\tilde f^m(x)$ and 
$\tilde f^m(x)\in\im\delta$ for all $m\ge 1$. Note that  
$\deg f>\deg \tilde f\ge 1$ (since $\deg f\ge 1$). 
Therefore, replacing $f$ by $\tilde f$ and repeating the same procedure, 
if necessary, we may assume $f'(x)\ne 0$. 
Consequently, by Eq.\,(\ref{FrobeniCase-peq1}) we  also have 
$\deg h_m\ge 1$ for all $m\ge 1$.

Now assume $p>2$. Then $f', \, (f^2)' \ne 0$, and 
by Eq.\,(\ref{FrobeniCase-peq1}) with $m=1, 2$  
we have $h_1',\, h_2'\ne 0$ and 
\begin{align}
h_2(x)-h_2^p(x)=(h_1(x)-h_1^p(x))^2=h_1^2(x)-2h_1^{p+1}(x)+h_1^{2p}(x). \label{FrobeniCase-peq2}
\end{align}
Applying $d/d x$ to the equation above we get 
\begin{align}
h_2'(x)= 2h_1(x)h_1'(x)-2h_1^p(x)h_1'(x).\label{FrobeniCase-peq3}
\end{align}

On the one hand, by the two equations above we have  
\begin{align}
\deg h_2&=2\deg h_1,\\
\deg h'_2&=p\deg h_1 +\deg h_1',
\end{align}
which imply 
$$
p\deg h_1 \le  p\deg h_1 +\deg h_1' =\deg h'_2 \le \deg h_2-1= 2\deg h_1-1.
$$ 
Since $ \deg h_1 \ge 1$, we get $p<2$, which is a contradiction. 

Therefore, we have $p=2$.  But in this case $f', \, (f^3)' \ne 0$, and 
by Eq.\,(\ref{FrobeniCase-peq1}) with $m=1, 3$  
we have $h_1',\, h_3'\ne 0$ and  
\begin{align}
2\deg h_3=3\deg f,
\end{align}
Applying $d/d x$ to Eq.\,(\ref{FrobeniCase-peq1}) with $m=3$ and $p=2$ we get 
\begin{align}
\deg h_3'(x)&= 2\deg f+\deg f',  
\end{align}
By the two equations above we  have $\deg f=2/3\deg h_3$ and 
$\deg h_3' = \frac43\deg h_3+\deg f'>\deg h_3$, 
which is a contradiction again. Therefore $\rad\big(\im \delta\big)=\{0\}$
\epfv

Next, the following example shows that 
the base field in the LFED Conjecture 
%#$ what about the LNED conjecture???  
can not be replaced by an integral domain 
of characteristic zero.

\begin{exam}\label{BaseDomainEx}
Let $t, x, y$ be commutative free variables; $R=\bC[t^{-1}, t]$; 
$\cA=R[x, y]$; and
$\phi\in\cEnd_R(\cA)$ that maps $x \to 2x$ and $y\to ty$. 
Then it is easy to verify that $\I-\phi$ is LF  
and $\im(\I-\phi)$ is the $R$-subspace spanned by $(1-2^a t^b)x^a y^b$ 
for all $a, b\in \bZ$. In particular, $x^m\in \im(\I-\phi)$ for all $m\ge 1$, 
since $(1-2^m)$ is invertible in $R$. But for all $m\ge 1$, 
$x^my \not \in \im(\I-\phi)$, since $(1-2^mt)$ is not invertible in $R$. 
Therefore, $\im(\I-\phi)$ is not a MS of $\cA$. 
\end{exam}
 
On the other hand, the following example shows that 
Problem \ref{LFNED-Prob} is also interesting for 
some algebras over an integral domain. 

\begin{exam}\label{phia-exam}
Let $a\in \bZ$, $x$ be a free variable, and 
$\phi_a$ the $\bZ$-algebra endomorphism of $\bZ[x]$ 
that maps $x$ to $ax$. Then $\im(\I-\phi_a)$ is the $\bZ$-subspace 
spanned by $(1-a^n )x^n$ for all $n\ge 1$. More precisely, 
\begin{align}\label{phia-exam-eq1}
\im(\I-\phi_a)=
\begin{cases}
x\bZ[x] &\text{ if } a=0;\\
\{0\}  &\text{ if } a=1;\\ 
2x\bZ[x^2] &\text{ if } a=-1; \\
\text{Span}_\bZ\{(1-a^n)x^n\,|\, n\ge 1\} &\text{ otherwise.}
\end{cases} 
\end{align}
The radical of $\im(\I-\phi)$ is given by 
\begin{align}\label{phia-exam-eq2}
\rad(\im(\I-\phi_a))=
\begin{cases}
x\bZ[x] &\text{ if } a=0;\\
\{0\}  &\text{ otherwise.}
\end{cases} 
\end{align}
Consequently, for all $a\in \bZ$ the image $\im(\I-\phi_a)$ is 
MS of $\bZ[x]$.
\end{exam}

\pf Eq.\,(\ref{phia-exam-eq1}) is obvious and the last statement can be easily verified by 
Eq.\,(\ref{phia-exam-eq2}) and Definition \ref{Def-MS}. 
To show Eq.\,(\ref{phia-exam-eq2}), the cases $a= 0, \pm 1$ are straightforward. 
So we assume $|a|\ge 2$. Note that $1-a^n$ in this case 
is not invertible in $\bZ$ for any $n\ge 1$. Note also that $\im(\I-\phi)$ is 
a homogeneous $\bZ$-subspace or $\bZ$-submodule of $\bZ[x]$.

Let $u\in\rad(\im(\I-\phi))$. Then $\deg u\ge 1$, for $\im(\I-\phi)$ 
obviously does not contain any nonzero constant. Replacing $u$ by a power of $u$ 
we assume $u^m\in \im(\I-\phi)$ for all $m\ge 1$. 
Let $bx^d$ be the leading term of $u$. Then  
$(bx^d)^m \in \im(\I-\phi)$ for all $m\ge 1$, whence  $(1-a^{md})\,|\, b^m$. 
Set $\tilde a\!:=a^d$. Then $(1-\tilde a^m)\,|\, b^m$ for all $m\ge 1$. 
Consequently, there are only finitely many distinct primes 
$p$ such that $p$ divides $1-\tilde a^m$ for some $m\ge 1$. 

On the other hand, for all co-prime $m, k\ge 1$ there 
exist $u(t), v(t)\in \bQ[t]$ such that 
$$
(1-t^m)u(t)+(1-t^k)v(t)=1-t.
$$ 
Furthermore, by going through the Euclidean algorithm for 
$1-t^m$ and $1-t^k$ it is easy to see that we can actually choose 
$u(t), v(t)\in \bZ[t]$.
 
Replacing $t$ by $\tilde a$ in the equation above we see that the integers 
$|(1-\tilde a^m)/(1-\tilde a)|$ and $|(1-\tilde a^k)/(1-\tilde a)|$ are co-prime for 
all co-prime $m, k\ge 1$. Since for all distinct $m\ge1$ the integers $1-\tilde a^m$ 
are all distinct, it is easy to see that 
there are infinitely many distinct primes $p$ such that $p$ 
divides $1-\tilde a^m$ for some $m\ge 1$. Contradiction.
\epfv

\renewcommand{\theequation}{\thesection.\arabic{equation}}
\renewcommand{\therema}{\thesection.\arabic{rema}}
\setcounter{equation}{0}
\setcounter{rema}{0}

\section{\bf Some Known Cases of the LFED and LNED Conjectures} \label{S3}

In this section we discuss some cases of conjectures \ref{LFED-Conj}, 
\ref{LNED-Conj} and \ref{Idem-Conj} that either are known in the literature  
or can be proved from some results in the literature.  
{\it Throughout this section $K$ denotes a field of characteristic zero
and $\cA$ a $K$-algebra.} 

We start with the following example. Although it is trivial,  
it can be read as a first positive sign for 
the LFED conjecture. 

\begin{exam}\label{TrivalEx}
Let $R$ be a unital commutative ring containing $\bQ$, $x$ a free variable and $D$ an 
arbitrary $R$-derivation of the polynomial algebra $R[x]$. Write 
$D=a(x)\frac{d}{d\, x}$ for some $a(x)\in R[x]$. Then 
$\im D$ is the principal ideal of $R[x]$ generated by $a(x)$,   
and hence a MS of $R[x]$.
\end{exam}

Furthermore, for the univariate polynomial algebra 
$K[x]$ the following theorem is 
proved  in \cite{OneVariableCase}.  

\begin{theo}\label{OneVarCase}
$1)$ The LFED conjecture  holds for all $K$-derivations   
and $K$-$\cE$-derivations (not necessarily LF) of $\kx$. 

$2)$ The LNED conjecture holds for 
all LN $K$-derivations of $\kx$.

$3)$ The LNED conjecture holds for 
all LN $K$-$\cE$-derivations $\delta$ of $\kx$ and the ideals $I$ of $\kx$ 
that are generated by a polynomial $u\in \kx$ with either $u=0$, or $\deg u\le 1$, 
or $u$ has at least one repeated root in the algebraic closure of $K$.   
\end{theo}

For multivariate polynomial algebras the following theorem is 
proved  in \cite{EWZ}, which can be 
re-stated as follows. 

\begin{theo}\cite[Theorem 3.1]{EWZ} \label{MultiVarCase}
The LFED conjecture holds for all LF $K$-derivations of 
the polynomial algebra over $K$ in two commutative free variables.
\end{theo}

For (multivariate) Laurent polynomial algebras the following theorem is 
proved in \cite{LaurentPolyCase}. 

\begin{theo}\label{MultiLaurentCase}
Let $x=(x_1, x_2, \dots, x_n)$ be $n$ commutative free variable and 
$K[x^{-1}, x]$ the Laurent polynomial algebra in $x$ over $K$. Then 
the following statements hold:
\begin{enumerate}
  \item[$1)$] $K[x^{-1}, x]$ has no nonzero locally nilpotent $K$-derivations or $K$-$\cE$-derivations. Hence, the LNED conjecture   holds (trivially) for $K[x^{-1}, x]$;
  \item[$2)$] if $n\le 2$, then the LFED conjecture   
              holds for $K[x^{-1}, x]$.
\end{enumerate}
\end{theo}

Next, we discuss some cases of the LFED and LNED conjectures 
for algebraic $K$-algebras. First, we have the following 

\begin{theo}\label{AlgLocCase}
Both the LFED conjecture  and 
the LNED conjecture hold for 
all local $K$-algebras $\cA$ 
that are algebraic over $K$.
\end{theo}

\pf Note first that by \cite[Theorem 7.6]{MS} the proper MSs 
of $\cA$ are characterized as follows: 
\vspace{2mm}
\begin{enumerate} 
\item[$(\ast)$] {\it a proper $K$-subspace $V$ of $\cA$ is a MS of $\cA$, 
if and only if $1\not \in V$}.  
\end{enumerate}
\vspace{2mm}
Let $\delta$ be a $K$-derivation or 
$K$-$\cE$-derivation of $\cA$.
Assume first that $\delta$ is LF.
If $1\not \in \im \delta$, then by 
the fact $(*)$ above $\im\delta$ is 
a MS of $\cA$. 
If $1\in \im \delta$, then $\im\delta=\cA$
by \cite[Proposition 1.4]{Idem}. 
Hence the LFED conjecture 
holds for $\cA$. 

Now assume that $\delta$ is LN. Then 
$1\in \cA^\delta$, if $\delta\in \Der_K(\cA)$, 
and by \cite[Proposition 2.1.3]{E} and \cite[Corollary 2.4]{Idem},  
it is also the case if $\delta\in \Eder_K(\cA)$. 
If $1 \in \im \delta$, then $\delta s=1$ for some $s\in \cA$. 
Then by \cite[Proposition 3.2]{Idem} $s$ is  
transcendental over $K$. Contradiction. 

Therefore $1\not \in \im \delta$. Then by 
the fact $(*)$ above $\delta$ actually maps every 
$K$-subspace of $\cA$ to a MS of $\cA$. In particular, the 
LNED conjecture holds for $\cA$ as well.
\epfv

Next, the following two theorems are proved  
in \cite{Alg-LFNED}. 
 
\begin{theo} \label{FinDimCase}
Both the LFED conjecture  and 
the LNED conjecture hold for 
all finite dimensional $K$-algebras.  
\end{theo}

\begin{theo}\label{AlgCase1}
Let $\cA$ be a $K$-algebra such that every $K$-subalgebra generated 
by finitely many elements of $\cA$ is finite dimensional, and  
$\delta$ a $K$-derivation or $K$-$\cE$-derivation of $\cA$.  
Then the following statements hold:
\begin{enumerate}
  \item[$1)$] if $\delta$ is LN, then $\delta$ maps every $K$-subspace 
of $\cA$ to a MS of $\cA$. In particular, the LNED conjecture holds for $\cA$; 
  \item[$2)$] if $\delta$ is a LF $K$-derivation, or a LF $K$-$\cE$-derivation of the form 
  $\delta=\I-\phi$ for some surjective $\phi\in\cEnd_K(\cA)$, 
  then $\im \delta$ is a MS of $\cA$, i.e.,  
the LFED conjecture holds for $\delta$.
\end{enumerate}
\end{theo}

For commutative algebraic $K$-algebras we here give a different proof for 
the proposition below, which is stronger for the $K$-derivation case 
than that of the theorem above. 

\begin{propo}\label{AlgCase2}
Let $\cA$ be a commutative $K$-algebra that is algebraic over $K$, and 
$\delta$ an arbitrary $K$-derivation, or a LN 
$K$-$\cE$-derivation of $\cA$. 
Then $\im\delta$ does not contain any nonzero idempotent of $\cA$. 
Consequently, $\delta$ maps every $K$-subspace of $\cA$ to a MS of $\cA$.
\end{propo}

\pf Let $e\in \im\delta$ be an idempotent, $u\in \cA$ 
such that $\delta u=e$, and $p(t)=t^d+\sum_{i=0}^{d-1} c_it^i \in K[t]$ 
the minimal polynomial of $u$.

We first consider the $K$-derivation case. 
Let $\delta=D\in \Der_K(\cA)$. 
Since $\cA$ is commutative, $De=D e^2=2eDe$, whence 
$(1-2e)De=0$. Note that $1-2e$ is a unit of $\cA$, 
for $(1-2e)^2=1$. Hence $De=0$. Consequently, 
$D^k(u^k)=k!\,e$ and $D^m(u^k)=0$ for all $m>k\ge 1$.  
Therefore $0=D^d p(u)=d!e$, whence $e=0$. 
Therefore $\im D$ 
does not contain any nonzero idempotent of $\cA$, and 
by Theorem \ref{IdemCri}, $D$ maps 
every $K$-subspace of $\cA$ to a MS of $\cA$. 

To show the $K$-$\cE$-derivation case, by a similar argument as above  
it suffices to show that $\delta e=0$ and $\delta^k (u^k)=k!\,e$ 
for all $k\ge 1$.

First, by van den Essen's one-to-one correspondence 
between the set of LN $K$-derivations of $\cA$ and the set of LN $K$-$\cE$-derivations 
of $\cA$ (see \cite{E0} or \cite[Proposition 2.1.3]{E}), there exists a LN 
$K$-derivation $D$ of $\cA$ such that $\delta=\I-e^D=\sum_{i=1}^\infty \frac{D^i}{i!}.$ 
By the argument in the first part of the second paragraph of 
this proof we have $De=0$, regardless $e\in \im D$ or not. 
Hence we also have $\delta e=0$.  

Next, we use the induction to show that $\delta^k(u^k)=k!e$ 
for all $k\ge 1$.

The case $k=1$ is trivial. Assume that 
$\delta^i(u^i)=i!e$ for all $1\le i\le k$. 
It is easy to check inductively that 
for all $x, y\in \cA$ and $n\ge 1$, we have 
\begin{align}\label{dlt-xy}
\delta^n(xy)=\sum_{i=0}^n \binom ni \delta^i(x)\,  (\I-\delta)^i
\big(\delta^{n-i}(y)\big).
\end{align}
Letting $n=k+1$, $x=u$ and $y=u^k$ in the equation above we have 
 \begin{align*}
 \delta^{k+1}(u^{k+1})&=\sum_{i=0}^{k+1} \binom {k+1}i \delta^i(u)\,  (\I-\delta)^i \big(\delta^{k-i+1}(u^k)\big).
\end{align*}
Then by the facts $\delta u=e$ and $\delta e=0$, the only nonzero term in the sum above 
is the $i=1$ term. Therefore by the induction assumption we have  
\begin{align*}
\delta^{k+1}(u^{k+1})=(k+1)e(\I-\delta)\big(\delta^k (u^k)\big)=(k+1)!e(\I-\delta)(e)=(k+1)!e.
\end{align*}
Hence, by induction we have $\delta^k(u^k)=k!e$ 
for all $k\ge 1$, as desired. 
\epfv

Note that the arguments in the proof of Proposition \ref{AlgCase2} above 
go through equally well for every $K$-derivation or $K$-$\cE$-derivation 
$\delta$ of $\cA$ and all idempotents $e \in \cA^\delta\cap \im\delta$, 
regardless of the commutativity of $\cA$. Therefore, we also have the following 

\begin{corol}
Let $\cA$ be a $K$-algebra (not necessarily commutative) that is algebraic over $K$, 
and $\delta$ an arbitrary $K$-derivation or $K$-$\cE$-derivation of $\cA$. 
Then $\cA^\delta \cap \im\delta$ does not contain any nonzero idempotent of $\cA$.
\end{corol}

For the $K$-algebras that are not algebraic over $K$, we have 
the following theorem proved in \cite{Idem}.

\begin{theo} %\label{MainThm-1}
Let $K$ be a field of characteristic zero and $\cA$ a $K$-algebra 
(not necessarily unital or commutative). Then the following statements hold:
\begin{enumerate}
  \item[$1)$] for every LF $D\in \Der_K(\cA)$ and an idempotent 
$e\in \cA^D\cap\im D$, we have $(e)\subseteq\im D$;
\item[$2)$]  for every LF $\delta\in \Eder_K(\cA)$ and an idempotent 
$e\in \cA^\delta\cap\im\delta$, we have $e\cA,\,\cA e  \subseteq\im\delta$. Furthermore, 
if $\delta$ is LN, we also have $(e) \subseteq\im\delta$. 
\end{enumerate}
\end{theo}

Note that, if $\cA$ is commutative, then for an arbitrary $K$-derivation or 
a LN $K$-$\cE$-derivation $\delta$ of $\cA$, we have that 
all idempotents of $\cA$ lie in $\cA^\delta$, as shown in the proof of 
Proposition \ref{AlgCase2}.  Therefore, from this fact and  
the theorem above we immediately 
have the following 

\begin{corol}%\label{Corol-1.1}
Assume that $\cA$ is commutative  (but not necessarily algebraic over $K$). 
Then Conjecture \ref{Idem-Conj} holds for all LF 
$D\in \Der_K(\cA)$ and all LN $\delta\in \cE$-$\Der_K(\cA)$. 
\end{corol}

The next case of Conjecture \ref{Idem-Conj} 
follows (somewhat unexpectedly) from the following 
classical {\it Singer-Wermer Theorem} in 
the theory of Banach algebras.  

\begin{theo}\label{D-BanachCase}
Let $\cA$ be a commutative Banach algebra  
and $D$ an arbitrary derivation. Then   $\im D$ 
is contained in the Jacobson radical $\mathfrak J(\cA)$ of $\cA$.   
\end{theo}

The theorem above was first proved by I. M. Singer and J. Wermer  
\cite{SW} in $1955$ for all continuous derivations, 
and in the same paper they also conjectured that the continuous condition 
is not necessary. More than thirty years later 
it was shown by M. P. Thomas \cite{T} in $1988$  
that it is indeed the case. 

Note that for all unital rings $R$ and nonzero 
idempotents $e\in R$, $1_R-e$ is also an idempotent, 
and can not be invertible. Then by \cite[Proposition 4.3]{Pi} the Jacobson 
radical $\mathfrak J(R)$ of $R$ does not contain any nonzero 
idempotent of $R$. From this general fact and 
Theorem \ref{D-BanachCase} we immediately have the following

\begin{corol}\label{D-BanachCaseCorol}
Let $\cA$ be a commutative Banach algebra  
and $D$ an arbitrary derivation. Then 
$\im D$ does not contain any nonzero idempotent $e$ of $\cA$. 
In particular, Conjecture \ref{Idem-Conj} holds for 
all derivations of $\cA$.  
\end{corol}
 
Note that there are also many results in the literature on the generalizations  
of the Singer-Wermer Theorem to certain derivations of some other  
algebras (e.g., see the survey paper \cite{MMa} and the book \cite[Section 6.4]{Pa}, 
and also the references therein). For example, it was shown in \cite{MR2} that every 
{\it centralizing} derivation $D$ (i.e., for all $u\in\cA$, $[u, Du]$ lies in 
the center of $\cA$) of an arbitrary Banach algebra $\cA$ has its image contained in the Jacobson radical of $\cA$. Hence Corollary \ref{D-BanachCaseCorol} and Conjecture \ref{Idem-Conj} 
also hold for centralizing derivations of all Banach algebras.

Conversely, Conjecture \ref{Idem-Conj} and more generally, 
the LFED and LNED conjectures in some sense provide 
some generalizations of the Singer-Wermer Theorem 
to (noncommutative) Banach algebras and also some other 
more general algebras (e.g., the normed algebras, etc.) 
in the general theory of Banach algebras.

\renewcommand{\theequation}{\thesection.\arabic{equation}}
\renewcommand{\therema}{\thesection.\arabic{rema}}
\setcounter{equation}{0}
\setcounter{rema}{0}

\section{\bf The LFED Conjecture from 
a Different Point of View} \label{S4}

Throughout this section {\it $K$ stands for a field of characteristic zero and $\cA$ for a $K$-algebra.}
In this section we discuss the LFED conjecture \ref{LFED-Conj} 
in terms of the decompositions of $\cA$ associated with 
the Jordan-Chevalley decompositions of LF $K$-derivations 
and $K$-$\cE$-derivations of $\cA$.

We first assume that $K$ is algebraically closed. 
For each LF $K$-linear endomorphism 
$\psi$ of $\cA$, let $\Lambda$ be the 
set of eigenvalues of $\psi$ and 
$\cA_\lambda\!:=\sum_{i=1}^\infty \Ker (\lambda\I-\psi)^i$ 
for all $\lambda \in \Lambda$.  Then it is well-known 
(e.g., see \cite[Proposition $1.3.8$]{E}),  
 \cite[Proposition $4.2$]{H}) that 
$\cA$ can be decomposed as 
\begin{align}\label{DeltaDecomp}
\cA=\oplus_{\lambda\in \Lambda}\cA_\lambda.
\end{align}

Furthermore, $\psi$ is said to be {\it semi-simple} if 
$\cA_\lambda$ $(\lambda\in \Lambda)$ coincides with the eigenspace of $\psi$ 
corresponding to the eigenvalue $\lambda$. 

With the decomposition as in Eq.\,(\ref{DeltaDecomp}) it can be readily 
verified (e.g., see the proof of \cite[Lemma 3.5 or 4.1]{Idem}) 
that the image $\im\psi$ can be decomposed as
\begin{align}\label{ImageDecomp}
\im\psi=\psi(\cA_0)\oplus \bigoplus_{0\ne \lambda\in \Lambda}\cA_\lambda.
\end{align}
 
If $\psi$ is a (LF) $K$-derivation of $\cA$, 
then by setting $\cA_\gamma=0$ for all $\gamma\not\in \Lambda$ we have $\cA_\lambda\cA_\mu\subseteq\cA_{\lambda+\mu}$ for all $\lambda, \mu\in \Lambda$, 
i.e., the decomposition in Eq.\,(\ref{DeltaDecomp}) is a so-called  
{\it additive algebra grading} of $\cA$. In particular, $\cA_0$ is 
a $\psi$-invariant $K$-subalgebra of $\cA$, and the restriction $\psi\,|_{\cA_0}$ 
is a LN $K$-derivation of the $K$-algebra $\cA_0$. Then by Eq.\,(\ref{ImageDecomp}) 
the image of the LF $K$-derivation $\psi$ of the $K$-algebra $\cA$ is completely 
determined by the image of the LN $K$-derivation $\psi\,|_{\cA_0}$ 
of the $K$-algebra $\cA_0$.
 
Similarly, if $\psi$ is a (LF) $K$-algebra endomorphism of $\cA$, 
then by setting $\cA_\gamma=0$ for all $\gamma\not\in \Lambda$ we have 
$\cA_\lambda\cA_\mu\subseteq\cA_{\lambda\mu}$ for all 
$\lambda, \mu\in \Lambda$, i.e., 
the decomposition in Eq.\,(\ref{DeltaDecomp}) is a so-called  
{\it multiplicative algebra grading} of $\cA$. 
In particular, $\cA_1$ is a $\psi$-invariant 
$K$-subalgebra of $\cA$, and the restriction $\psi\,|_{\cA_1}$ 
is a $K$-algebra endomorphism of the $K$-algebra $\cA_1$ 
such that $\I_{\cA_1}-\psi\,|_{\cA_1}$ is a LN 
$K$-$\cE$-derivation of $\cA_1$. 
 
Now set $\delta\!:=\I-\psi$. Then $\delta$ is a LF $K$-$\cE$-derivation of $\cA$, 
and $\delta\,|_{\cA_1}$ is a LN $K$-$\cE$-derivation of $\cA_1$.  
By Eq.\,(\ref{ImageDecomp}) with $\psi$ replaced by $\delta$ we have  
\begin{align}\label{DeltaImageDecomp}
\im\delta=\delta(\cA_1)\oplus \bigoplus_{1\ne \lambda\in \Lambda} \cA_\lambda.
\end{align}
Therefore, the image of the LF $K$-$\cE$-derivation $\delta$ of the $K$-algebra $\cA$ 
is completely determined by the image of the LN $K$-$\cE$-derivation $\delta\,|_{\cA_1}$ 
of the $K$-algebra $\cA_1$.

One special but important case is when the $K$-linear map $\psi$ is a semi-simple $K$-derivation 
or a semi-simple $K$-endomorphism of $\cA$, i.e., $\cA_\lambda$ $(\lambda\in \Lambda)$ 
in Eq.\,(\ref{DeltaDecomp}) is the eigenspace of $\psi$ corresponding to the eigenvalue 
$\lambda$ of $\psi$. In this case the $K$-subspaces $\psi(\cA_0)$ in Eq.\,(\ref{ImageDecomp}) 
and $\delta(\cA_1)$ in Eq.\,(\ref{DeltaImageDecomp}) are respectively equal to zero. 
Based on this observation and also the LFED conjecture we propose the following 
what we call {\it the Grading conjecture}. 

\begin{conj}$(${\bf The Grading Conjecture}$)$ \label{Grading-Conj}
Let $K$ be a field of characteristic zero (not necessarily algebraically closed), 
$\cA$ a $K$-algebra and $(\Lambda, \ast)$ a monoid  
with the unit $e$. Assume that $\cA$ has the following $K$-algebra 
grading with respect to the monoid $(\Lambda, \ast)$:  
\begin{align}\label{DeltaDecomp2}
\cA=\oplus_{\lambda\in \Lambda}\cA_\lambda,
\end{align}
i.e., $\cA_\lambda$ is a $K$-subspace of $\cA$ for each $\lambda\in \Lambda$,  
and $\cA_\lambda\cA_\mu\subseteq \cA_{\lambda\ast\mu}$ for all 
$\lambda, \mu\in \Lambda$. Then the $K$-subspace 
$\bigoplus_{e \ne \lambda\in \Lambda} \cA_\lambda$ is a MS of $\cA$. 
\end{conj}

Two remarks on the conjecture above are as follows.

First, it is easy to see that the set of all semi-simple $K$-derivations (resp., $K$-algebra endomorphisms) of $\cA$ is in one-to-one correspondence with the set of all additive (resp., multiplicative) algebra gradings of $\cA$ with the index monoid $(\Lambda, \ast)$ being a sub-monoid of $(K, +)$ (resp. $(K, \cdot)$). 
Therefore, the semi-simple $K$-derivation (resp., $K$-$\cE$-derivation) case of the LFED conjecture is equivalent to the case of the Grading conjecture under the extra assumption 
that the index monoid $(\Lambda, \ast)$ is a sub-monoid of $(K, +)$ (resp. $(K, \cdot)$).        
For convenience, we refer to these two cases of the Grading conjecture respectively 
as {\it the additive Grading conjecture} and 
{\it the multiplicative Grading conjecture}.

Second, by the equivalences mentioned above, some of the known cases of the LFED conjecture 
discussed in Section \ref{S3} and also some of those that will be 
proved in the next two sections can be translated or re-formulated as 
certain cases of the additive and multiplicative Grading conjectures. 
For example, both the additive and multiplicative Grading conjectures hold 
for the univariate polynomial algebra $K[x]$ by Theorem \ref{OneVarCase};  
and all local algebraic $K$-algebra by Theorem \ref{AlgLocCase}; and 
all finite dimensional $K$-algebras by Theorem \ref{FinDimCase}; etc..
They also hold for the $K$-algebras $\cA$ in Theorem \ref{AlgCase1}, 
which can be shown in the following 

\begin{corol} 
Let $\cA$ be as in Theorem \ref{AlgCase1}. Assume that $\cA$ has a $K$-linear 
decomposition as in Eq.\,(\ref{DeltaDecomp}) (with $\Lambda\subseteq K$). 
Then the following statements hold:
\begin{enumerate}
  \item[$1)$] if the decomposition in Eq.\,(\ref{DeltaDecomp}) is an additive 
   algebra grading of $\cA$, then  
the $K$-subspace $\bigoplus_{0\ne \lambda\in \Lambda} \cA_\lambda$
is a MS of $\cA$;  
  \item[$2)$] if the decomposition in Eq.\,(\ref{DeltaDecomp}) is a 
  multiplicative algebra grading of $\cA$, then    
the $K$-subspace $\bigoplus_{1\ne \lambda\in \Lambda} \cA_\lambda$
is a  MS of $\cA$.    
\end{enumerate} 
In other words, both the additive and multiplicative Grading conjectures 
hold for $\cA$. 
\end{corol}

\pf $1)$ Define $D:\cA\to\cA$ by setting $Du=\lambda u$ for all $\lambda\in\Lambda$ and 
$u\in\cA_\lambda$. Since the decomposition in Eq.\,(\ref{DeltaDecomp}) is an additive 
algebra grading of $\cA$, it is easy to see that $D$ is a LF $K$-derivation of $\cA$ 
with  $\im D=\bigoplus_{0\ne \lambda\in \Lambda} \cA_\lambda$. 
Then by Theorem \ref{AlgCase1}, $2)$ the statement follows. 

$2)$ Define $\phi:\cA\to\cA$ by setting $\phi (u)=\lambda u$ for all $\lambda\in\Lambda$ and 
$u\in\cA_\lambda$. Since the decomposition in Eq.\,(\ref{DeltaDecomp}) is a multiplicative algebra grading of $\cA$, it is easy to see that $\phi$ is a LF $K$-algebra endomorphism 
of $\cA$. Set $\delta\!:=\I-\phi$. Then $\delta$ is a LF $K$-$\cE$-derivation of $\cA$ with 
$\im \delta =\bigoplus_{1 \ne \lambda\in \Lambda} \cA_\lambda$. 

Note that $\cA_0=\Ker\phi$ and is an ideal of $\cA$. Set $\bar\cA\!:=\cA/\cA_0$ and 
$\bar\delta\!:=\I_{\bar\cA}-\bar\phi$, where $\bar\phi$ is the $K$-algebra 
endomorphism of $\bar\cA$ induced by $\phi$.  
We may identify $\bar\cA$ with the $K$-subalgebra 
$\bigoplus_{0 \ne \lambda\in \Lambda} \cA_\lambda$ of $\cA$. 
Then under this identification $\bar\phi(u)=\lambda u$ for all 
$0\ne \lambda\in\Lambda$ and $u\in\cA_\lambda$. In particular, 
$\bar\phi$ is a LF $K$-algebra automorphism 
of $\bar\cA$ and $\bar\delta=\I_{\bar\cA}-\bar\phi$ is a LF 
$K$-$\cE$-derivation of $\bar\cA$. Then by 
Theorem \ref{AlgCase1}, $2)$ $\im\bar \delta$ is 
a MS of $\bar\cA$. Note that $\im\bar \delta=\im \delta/\cA_0$ and the ideal 
$\cA_0$ is obviously contained in $\im \delta$. Then   
by \cite[Proposition 2.7]{MS} $\im\delta$ 
is a MS of $\cA$, whence the statement follows. 
\epfv

Besides the cases above, we also have 
the following cases of the additive and multiplicative 
Grading conjectures in a more general setting.    
 
\begin{propo}
Let $K$ be a field of characteristic zero  
and $\cA$ a $K$-algebra with a decomposition 
as in Eq.\,(\ref{DeltaDecomp}) (with $\Lambda\subseteq K$). 
Let $H_1$ (resp., $H_2$) be the sub-monoid of 
the abelian group $(K, +)$ (resp., $(K\backslash\{0\}, \cdot)$) 
generated by elements $0\ne \lambda\in \Lambda$ (resp., $0, 1\ne \lambda\in\Lambda$).  
Then the following statements hold:
\begin{enumerate}
  \item[$1)$] if the decomposition in Eq.\,(\ref{DeltaDecomp}) is an additive 
   algebra grading of $\cA$ and $0\not\in H_1$, then for every $\vartheta$-MS\, $V$ of $\cA_0$,  
the $K$-subspace $V \oplus\bigoplus_{0\ne \lambda\in \Lambda} \cA_\lambda$
is a $\vartheta$-MS of $\cA$;  
  \item[$2)$] if the decomposition in Eq.\,(\ref{DeltaDecomp}) is a 
  multiplicative algebra grading of $\cA$  and $1\not\in H_2$, then for every $\vartheta$-MS\, $V$ of $\cA_1$, the $K$-subspace $V \oplus\bigoplus_{1\ne \lambda\in \Lambda} \cA_\lambda$
is a $\vartheta$-MS of $\cA$.    
\end{enumerate}
\end{propo}

\pf Note that the $K$-subspace $\bigoplus_{0\ne \lambda\in \Lambda} \cA_\lambda$ in statement $1)$ under the condition $0\not\in H_1$ is an ideal of $\cA$, and the same for the $K$-subspace $\bigoplus_{1\ne \lambda\in \Lambda} \cA_\lambda$ in statement $2)$ under the condition $1\not\in H_2$. Then both statements $1)$ and $2)$ follow directly from \cite[Proposition 2.7]{MS}.  
\epfv

Next, we discuss an important special case of the multiplicative Grading conjecture. 
Let $z=(z_1, z_2, \dots, z_n)$ be $n$ commutative or noncommutative free variables and 
$\cA[z^{-1}, z]$ the algebra of Laurent polynomials  
in $z$ over a $K$-algebra $\cA$. Let $q=(q_1, q_2, \dots, q_n) \in K^n$ be 
such that $0\ne q_i\in K$ $(1\le i\le n)$ and $q^\alpha\!:=\prod_{i=1}^n q_i^{\alpha_i}\ne 1$ for 
all $0\ne \alpha=(\alpha_1, \alpha_2, \dots, \alpha_n)\in \bZ^n$ (e.g., let 
$q_i$ $(1\le i\le n)$ be $n$ distinct prime integers).  

Let $\Lambda=\{q^\alpha\,|\, \alpha\in \bZ^n\}$ and 
$V_\lambda$ ($\lambda \in \Lambda$) be the $K$-subspace 
formed by all $f(z)\in \cA[z^{-1}, z]$ such that $f(q_1z_1,\, q_2z_2, \dots, \, q_nz_n)=\lambda f(z)$.
Then it is easy to see that $\cA[z^{-1}, z]$ can be decomposed as 
\begin{align}\label{q-Decom-Laurent}
\cA[z^{-1}, z]=\oplus_{\lambda\in \Lambda} V_\lambda,
\end{align}
which is a multiplicative algebra grading of $\cA[z^{-1}, z]$ with $V_0=\cA$. 

Then the multiplicative Grading conjecture for 
the multiplicative $K$-algebra grading of $\cA[z^{-1}, z]$ in Eq.\,(\ref{q-Decom-Laurent}) 
becomes the following  

\begin{conj}\label{GeneralDK-Conj}
Let $K$ be a field of characteristic zero, $\cA$ a $K$-algebra, and $z=(z_1, z_2, \dots, z_n)$ $n$ 
commutative or noncommutative free variables. 
Denote by $M$ the $K$-subspace of the Laurent polynomial algebra 
$\cA[z^{-1}, z]$ consisting of the Laurent polynomials with no constant term.  
Then $M$ is a MS of $\cA[z^{-1}, z]$. 
\end{conj}

One known case of the conjecture above is as follows. 
Let $z_i$ $(1\le i\le n)$ be commutative free variables,  
$\cA=K$. Then the conjecture above in this case 
follows directly from the following remarkable Duistermaat-van der Kallen 
Theorem \cite{DK}, which is also the special case of 
the Mathieu conjecture \cite{Ma} for complex tori.
 
\begin{theo}\label{DvK-thm}
Let $K$ be a field of characteristic zero, $z=(z_1, \dots, z_n)$ 
commutative free variables and $M$ the $K$-subspace of $K[z^{-1}, z]$ 
of the Laurent polynomials with no constant term. Then 
$\rad(M)$ consists of $f\in  K[z^{-1}, z]$ such that 
$0$ does not lie in the polytope of $f$. 
Consequently, $M$ is a MS of $K[z^{-1}, z]$.  
\end{theo}

From the discussion above, we see that the conjecture \ref{GeneralDK-Conj} 
and more generally, the Grading conjecture \ref{Grading-Conj} can be viewed as 
some natural generalizations of the Duistermaat-van der Kallen 
Theorem.

\renewcommand{\theequation}{\thesection.\arabic{equation}}
\renewcommand{\therema}{\thesection.\arabic{rema}}
\setcounter{equation}{0}
\setcounter{rema}{0}

\section{\bf The LFED Conjecture for Some Special $\cE$-Derivations} \label{S5}

Throughout this section {\it $R$ denotes a unital commutative ring 
and $\cA$ an $R$-algebra}. We denote by $\nil(\cA)$ the 
set of all nilpotent elements of $\cA$ 
(although $\cA$ may not be commutative).

We shall show the LFED Conjecture \ref{LFED-Conj} 
for the $R$-$\cE$-derivations associated with some special 
$R$-algebra endomorphisms of $\cA$. We start with the following lemma, 
which will also play an important role in the next section of this paper.

\begin{lemma}\label{GenralLma}
Let  $A, B, C, D$ be four commutating $R$-module endomorphisms of 
$\cA$ such that $AB=0$ and $AD+BC=\I$. Then $\im A=\Ker B$.  
\end{lemma}

\pf Since $BA=AB=0$, we have $\im A\subseteq \Ker B$. Now let $a\in \Ker B$. 
Then $a=(AD+BC)(a)=(AD+CB)(a)=A(D(a))$. Therefore $a\in \im A$, 
whence the lemma follows. 
\epfv

Next, we consider $R$-$\cE$-derivations associated with 
{\it $R$-projections} (i.e., $\phi\in \cEnd_R(\cA)$ with $\phi^2=\phi$) 
and {\it $R$-involutions} (i.e., $\phi\in \cEnd_R(\cA)$ with  
$\phi^2=\I$) of $\cA$.

\begin{propo}\label{ProjInvoCases}
Let $\phi$ be an $R$-algebra endomorphism of $\cA$. Then the following statements hold: 
\begin{enumerate}
\item[$1)$] if $\phi^2=\phi$, then $\im(I-\phi) = \Ker \phi$.
\item[$2)$] if $\phi^2 = \I$ and $2\cdot 1_R$ is a unit of $R$, 
then  
\begin{align}
\im(\I-\phi)&=\Ker (\I+\phi),\label{ProjInvoCases-eq1}\\
\rad(\im(\I-\phi))&=\nil (\cA). \label{ProjInvoCases-eq2}
\end{align}
%where $\nil(\cA)$ denotes the set of all nilpotent elements of $\cA$.
\end{enumerate}

In both cases above, $\im(I-\phi)$ is a MS of $\cA$.
\end{propo}

\pf $1)$ Since $\phi^2 = \phi$, by Lemma \ref{GenralLma} above with 
$A=\I-\phi$, $B=\phi$ and $C=D=\I$, we have $\im(\I-\phi) = \Ker \phi$, 
which is an ideal of $\cA$, and hence also a MS of $\cA$.

$2)$ Since $2\cdot 1_R$ is a unit of $R$, we may apply Lemma \ref{GenralLma} 
with $A=\I-\phi$, $B=\I+\phi$ and $C=D=\frac12\I$, from which we get 
Eq.\,(\ref{ProjInvoCases-eq1}).

Now let $a \in \rad(\im(\I-\phi))$. Replacing $a$ by a power of $a$ 
we may assume that $a, a^2 \in \im(\I-\phi)$. Then by Eq.\,(\ref{ProjInvoCases-eq1})
we have $a, a^2 \in \Ker(\I+\phi)$, whence 
$\phi(a) = -a$ and
$\phi(a^2) = -a^2$. But since $\phi$ is an $R$-algebra endomorphism of $\cA$, 
we also have $\phi(a^2) = \phi(a)^2 = (-a)^2 = a^2$. Hence $a^2 = -a^2$ and $a^2 = 0$, 
for $2\cdot 1_R$ is a unit of $R$. Therefore $a\in \nil(\cA)$.  
Since $\nil (\cA)$ is obviously contained in $\rad(\im(\I-\phi))$, we have 
$\rad(\im(\I-\phi)) =\nil (\cA)$, and by Definition \ref{Def-MS},  
$\im(\I-\phi)$ is a MS of $\cA$. 
\epfv
 
Note that Proposition \ref{ProjInvoCases}, 
$2)$ is not always true if $2\cdot 1_R$ is not a unit of $R$ 
(e.g., see Example \ref{p-Count-LFED}).

Next, we show the following lemma, which reduces the $\cE$-derivation case 
of the LFED conjecture to 
the LF $\cE$-derivations associated with injective 
algebra endomorphisms.  

\begin{lemma}\label{ReducLma}
Let   $\phi$ be an 
$R$-algebra endomorphism of $\cA$. Set $\Ker_{\ge 1} \phi \!:=\sum_{i\ge 1} \Ker \phi^i$ 
%(and $\delta=\I-\phi$). 
and $\bar\cA\!:= \cA/\Ker_{\ge 1} \phi$. Denote by $\pi$ the quotient map from $\cA$ to $\bar\cA$ and $\bar\phi$ the induced map of $\phi$ 
from $\bar\cA$ to $\bar\cA$. Then 
\begin{enumerate}
  \item[$1)$] $\Ker_{\ge 1} \phi \subseteq \im(\I-\phi)$.
  \item[$2)$]  $\bar\phi$ is injective. 
\item[$3)$]  The following  equations hold:
\begin{align}\label{ReducLma-eq0}
\pi\big(\im(\I_\cA-\phi)\big)=\im (\I_{\bar\cA}-\bar\phi).
\end{align}
\begin{align}\label{ReducLma-eq1}
\im(\I_\cA-\phi)=\pi^{-1}\big(\im (\I_{\bar\cA}-\bar\phi)\big).
\end{align} 
\begin{align}\label{ReducLma-eq2}
\rad\big( (\im(\I_\cA-\phi) \big)=\pi^{-1}\big
(\rad\big(\im (\I_{\bar\cA}-\bar\phi)\big)\big).
\end{align} 
\item[$4)$] $ \im(\I-\phi)$ is a MS of $\cA$, if and only if 
$\im (\I_{\bar\cA}-\bar\phi)$ is a MS of $\bar\cA$.
  \end{enumerate}
\end{lemma}

\pf $1)$ Let $a\in \Ker_{\ge 1} \phi$. Then $\phi^k(a) = 0$ for some $k \ge 1$. 
Let $v=\sum_{i=0}^\infty \phi^i(a)$, which is a well-defined element of $\cA$. 
Then $(\I-\phi)v=(\I-\phi)(\sum_{i=1}^\infty \phi^i)(a)=a$. Therefore $a\in \im (\I-\phi)$. 
 
$2)$ Let $a\in \cA$ such that $\bar \phi(\pi(a))=0$. Since $\bar\phi\pi=\pi\phi$, 
we have $\pi(\phi(a))=0$, i.e., $\phi(a)\in \Ker_{\ge 1}\phi$. Then  
$\phi^{k+1}(a)=\phi^k(\phi(a))=0$ for some $k\ge 1$. Therefore,     
$a\in \Ker_{\ge 1}\phi=\Ker\pi$, whence $\pi(a)=0$ 
and $\bar\phi$ is injective.

$3)$ Since $\pi\phi=\bar\phi\pi$, we have 
$\pi(\I_\cA-\phi)=(\I_{\bar\cA}-\bar\phi)\pi$, from 
which and the surjectivity of $\pi$ we have Eq.\,(\ref{ReducLma-eq0}).
To show  Eq.\,(\ref{ReducLma-eq1}), first, by Eq.\,(\ref{ReducLma-eq0}) 
we have $\im(\I_\cA-\phi)\subseteq \pi^{-1}\big(\im (\I_{\bar\cA}-\bar\phi)\big)$. 
Let $a\in \pi^{-1}\big(\im (\I_{\bar\cA}-\bar\phi)\big)$. Then  
$\pi(a)\in \im (\I_{\bar\cA}-\bar\phi)$, i.e., there exists $b\in \cA$ 
such that 
$$
\pi(a)=(\I_{\bar\cA}-\bar\phi)(\pi(b))=\pi (\I_\cA-\phi)(b).
$$
Set $c\!:=(\I_\cA-\phi)(b)$. Then $c\in \im (\I_\cA-\phi)$ 
and $a-c\in \Ker \pi$. Since $\Ker\pi=\Ker_{\ge 1}\phi$, 
by statement $1)$ we have $a-c\in \im (\I_\cA-\phi)$. 
Hence $a=(a-c)+c\in  \im (\I_\cA-\phi)$, and Eq.\,(\ref{ReducLma-eq1}) follows.

To show Eq.\,(\ref{ReducLma-eq2}), 
first by Eq.\,(\ref{ReducLma-eq0}) we immediately have 
\begin{align*}
\pi \big(\rad \big( (\im(\I_\cA-\phi) \big)\big) \subseteq  
 \rad\big(\im (\I_{\bar\cA}-\bar\phi)\big).\\ 
\rad \big( (\im(\I_\cA-\phi) \big) \subseteq \pi^{-1}\big
(\rad\big(\im (\I_{\bar\cA}-\bar\phi)\big)\big). 
\end{align*}

Now let $a\in \pi^{-1}\big
(\rad\big(\im (\I_{\bar\cA}-\bar\phi)\big)\big)$. Then 
$\bar a\!:=\pi(a)\in \rad\big(\im (\I_{\bar\cA}-\bar\phi)\big)$, i.e., 
$\pi(a^m)=\bar a^m\in \im (\I_{\bar\cA}-\bar\phi)$, 
and hence $a^m\in \pi^{-1}\big(\im (\I_{\bar\cA}-\bar\phi)\big)$, 
for all $m\gg 0$.  
Then  by Eq.\,(\ref{ReducLma-eq1}), $a^m\in \im(\I_\cA-\phi)$ for all $m\gg 0$. 
Hence $a\in \rad \big( (\im(\I_\cA-\phi) \big)$ and Eq.\,(\ref{ReducLma-eq2}) follows. 

$4)$ follows directly from statement $1)$, Eq.\,(\ref{ReducLma-eq0})  and 
Proposition $2.7$ in \cite{MS}.
\epfv

Now we consider the following special family of $\cE$-derivations.

\begin{propo}\label{phi-i=phi-j}
Assume that $\cA$ is commutative and torsion-free as a $\bZ$-module, 
i.e., no $0\ne m\in \bZ$ is a zero-divisor of $\cA$.  
Let $\phi\in \cEnd_R(\cA)$ such that $\phi^i=\phi^j$ 
for some $1\le i<j$. Set $\delta\!:=\I-\phi$ and 
$\Ker_{\ge 1}\phi\!:=\sum_{k\ge 1}\Ker \phi^k$. Then   
\begin{align}
\rad(\im \delta)=\rad(\Ker \phi^i)=\rad(\Ker_{\ge 1}\phi).
\end{align} 
Consequently, $\im \delta$ is a MS of $\cA$.
\end{propo}

\pf First, the case $\phi=0$ or $\I$ is trivial. 
So we assume $\phi\ne 0, \I$. Second, 
since $\phi^i=\phi^j$ with $i<j$, we have $\phi^i=\phi^m$ for all $m\ge 1$ 
of the form $m=i+q(j-i)$ $(q\ge 0)$. 
Then for each $k\ge i$, choosing $q$ large enough such that $k\le m\!:=i+q(j-i)$ 
we have
$$
\Ker \phi^i\subseteq \Ker\phi^k\subseteq \Ker\phi^m=
\Ker\phi^i.
$$
Hence $\Ker\phi^i= \Ker\phi^k$ for all $k\ge i$ and 
$\Ker\phi^i=\Ker_{\ge 1}\phi$.  

Let $\pi$ be the quotient map from $\cA$ to 
$\bar \cA\!:=\cA/\Ker_{\ge 1}\phi$, and 
$\bar\phi$ the $R$-algebra endomorphism of $\bar \cA$ 
induced by $\phi$.  
Since $\pi^{-1}(\nil(\bar \cA))=\rad(\Ker\pi)
=\rad(\Ker_{\ge 1}\phi)$,  
by Eq.\,(\ref{ReducLma-eq2}) it suffices to show 
$\rad\big(\im(\I_{\bar\cA}-\bar\phi)\big)=\nil (\bar\cA)$.

Furthermore, by replacing $\cA$ by $\bar\cA$ and $\phi$ by $\bar\phi$, and 
by Lemma \ref{ReducLma}, $2)$ we may assume that $\phi$ is injective, and 
only need to show the following equation:
\begin{align}\label{OnlyEq}
\rad(\im\delta)=\nil (\cA). 
\end{align}
 
First, by Definition \ref{Def-rad} 
$\nil (\cA)$ is obviously contained in $\rad(\im\delta)$. 
Conversely, let $a\in \rad(\im \delta)$. Replacing $a$ by  
a power of $a$ we assume that $a^m\in \im\delta$ 
for all $m\ge 1$. 

Second, under the injective assumption on $\phi$,  
the condition $\phi^i=\phi^j$ with $i<j$ implies $\I=\phi^n$, 
where $n\!:=j-i$. Then $n\ge 2$, for we have assumed 
$\phi\ne \I$.

Since $\I-\phi^n=0$, we have $g(\phi)\delta=g(\phi)(\I-\phi)=0$, where 
$g(t)\!:=\sum_{k=0}^{n-1} t^i$. Hence $\im \delta \subseteq \Ker g(\phi)$.  
Therefore $a^m\in \Ker g(\phi)$ for all $m\ge 1$. 
Set $b_i=\phi^i(a)$ for all $0\le i\le n-1$. 
Then for all $m\ge 1$, we have
\begin{align}\label{phi-i=phi-j-peq3}
b_0^m +b_1^m+\cdots+b_{n-1}^m=0.
\end{align}

Note that the left-hand side of Eq.\,(\ref{phi-i=phi-j-peq3}) 
is the value at $b_i$ $(0\le i\le n-1)$ of the $m$-th power sum 
symmetric polynomial $p_m(x)\!:=\sum_{i=0}^{n-1}x_i^m$. 
It is well-known (e.g., see \cite{Mac}, or \cite{Wiki} and references therein.) that 
each elementary symmetric polynomial $e_m$ 
$(m\ge 1)$ can be written as a polynomial in $p_m$ $(m\ge 1)$ with coefficients 
in $\bQ$. Therefore, for all $m\ge 1$, the values $e_m(b_i; 0\le i\le n-1)$ (in $\cA$) of $e_m$  
at $b_i$ $(0\le i\le n-1)$, when viewed as elements of $\bQ\otimes_\bZ\cA$, 
are all equal to zero.  

On the other hand, since $\cA$ is a torsion-free $\bZ$-module
(and $\bZ$ is a PID), $\cA$ is also a flat $\bZ$-module 
(e.g., see \cite[Chapter I, $\S$2.4, Prop. 3]{Bou}). 
In particular, the homomorphism $\cA=\bZ\otimes_\bZ\cA \to \bQ\otimes_\bZ\cA$  
is injective. Therefore, $e_m(b_i; 0\le i\le n-1)$ for all $m\ge 1$  
(when viewed as elements of $\cA$) are also equal to zero.    
Consequently, $\prod_{i=0}^{n-1}(t-b_i)=t^n$ in $\cA[t]$. 
Letting $t=b_0=a$ we get $a^n=0$, whence $a\in \nil(\cA)$, 
as desired.
\epfv  

From Proposition \ref{phi-i=phi-j} or from its proof above 
we immediately have the following 

\begin{corol}\label{FinOrdCase}
Assume that $\cA$ is commutative and torsion-free as a $\bZ$-module. 
Then for every finite order $R$-algebra automorphism $\phi$ of $\cA$, we have 
\begin{align}\label{FinOrdCase-eq1}
\rad(\im (\I-\phi)) = \nil(\cA).
\end{align}
In particular, $\I-\phi$ maps every $R$-subspace 
of $\cA$ to a MS of $\cA$.
\end{corol}

\renewcommand{\theequation}{\thesection.\arabic{equation}}
\renewcommand{\therema}{\thesection.\arabic{rema}}
\setcounter{equation}{0}
\setcounter{rema}{0}

\section{\bf Some Cases for Algebraic Derivations and $\cE$-Derivations of Domains}\label{S6}

Throughout this section {\it $R$ stands for a unital commutative ring that 
contains $\bZ$ as a subring, and $\cA$ a unital $R$-algebra that is torsion-free 
as a $\bZ$-module. }
For convenience, we also assume $\bZ\subseteq R\subseteq \cA$. 
If $\cA$ has no left or right zero-divisors, 
we say $\cA$ is a {\it domain}.

Recall that an $R$-derivation or $R$-($\cE$-)derivation 
$\delta$ of $\cA$ is {\it algebraic over $R$}   
if there exists a nonzero polynomial $f(t) \in R[t]$ 
such that $f(\delta)=0$. When the base ring $R$ is clear in the context, 
we also simply say that $\delta$ is {\it algebraic}.

In this section we mainly consider some cases of Problem \ref{LFNED-Prob} 
for algebraic derivations and $\cE$-derivations of domains. 
In particular, we show that both the LFED conjecture and the LNED conjecture 
hold for all LF or LN algebraic derivations and $\cE$-derivations 
of integral domains of characteristic zero (see Theorem \ref{AlgIntDomainThm}). 
The proof will be divided into several lemmas and propositions, 
some of which will be proved  in more general settings.

\begin{lemma}\label{NoNilHD}
Assume further that $\cA$ is reduced, i.e., $\cA$    
has no nonzero nilpotent element.  
Then $\cA$ has no nonzero nilpotent $R$-derivations 
or $R$-$\cE$-derivations.
\end{lemma}

\pf Here, we only show the $R$-$\cE$-derivation case.  
The $R$-derivation case can be proved similarly. 

Assume otherwise and let $\phi\in \cEnd_R(\cA)$ such that the 
$R$-$\cE$-derivation $\delta\!:=\I-\phi$ is nonzero and nilpotent. 
Let $k\ge 2$ be the least positive integer such that $\delta^k=0$. 
Then there exists $u\in \cA$ such that $\delta^{k-1}u\ne 0$. 

By Eq.\,(\ref{dlt-xy}) it is easy to see 
that for all $m\ge 1$ and $v\in \cA$ with $\delta^2 v=0$, we have
\begin{align}
\delta^m\big( uv)  =(\delta^mu)v+m(\delta^{m-1}u-\delta^mu)\delta v.
\end{align}

Then by letting $m=k$ and $v=\delta^{k-2}u$, and applying the assumption 
$\delta^k=0$ we get   
\begin{align*}
0&=\delta^k\big( u(\delta^{k-2}u) \big)=k(\delta^{k-1}u)^2.
\end{align*}
Since $\cA$ is reduced and torsion-free as a $\bZ$-module,  
we have $\delta^{k-1}u=0$. Contradiction.
\epfv

Next, let us recall the following proposition proved 
in \cite[Theorem $4.6$]{RadKerDiff}. %\label{NoncommDomainCase-1} there.

\begin{propo}\label{NoAlgDer}  
Let $R$ be a unital integral domain of characteristic zero 
and $\cA$ a unital reduced $R$-algebra (not necessarily commutative) 
that is torsion-free as an $R$-module. 
Then $\cA$ has no nonzero $R$-derivation that is locally  
algebraic over $R$. 
In particular, $\cA$ has no nonzero $R$-derivation that is   
algebraic over $R$.  
\end{propo}

\begin{rmks}
$1)$ Proposition \ref{NoAlgDer} does not always hold for $\cE$-derivations, e.g., 
taking $\phi$ to be a non-identity finite order automorphism of $\cA$, 
if there is any. 

$2)$ Assume further that $\cA$ is a domain of characteristic zero. Then 
by Proposition \ref{NoAlgDer}, both the LFED and LNED Conjectures   
hold (trivially) for $R$-derivations of $\cA$ that 
are algebraic over $R$.   
\end{rmks}  

Next we consider algebraic $\cE$-derivations of domains of characteristic zero.

\begin{lemma}\label{(1-t)Lma}
Assume further that $R$ is an integral domain of characteristic zero, 
and  $\cA$ is a domain (containing $R$). 
Let $0,\I\ne \phi\in \cEnd_R(\cA)$ be algebraic over $R$, and 
$f(t)$ a minimal polynomial of $\phi$, i.e., $f(t)$ has the 
least degree among all $0\ne g(t)\in R[t]$ with $g(\phi)=0$.  
Then $f(t)=(1-t)h(t)$ for some 
$h(t)\in R[t]$ with $\deg h\ge 1$ and $h(1)\ne 0$.   
\end{lemma}

\pf Let $K_R$ be the field of fractions of $R$ and 
$\bar K_R$ be the algebraic closure of $K_R$. 
Decompose $f(t)$ in $\bar K_R[t]$ as   
\begin{align}\label{Decom-ft}
f(t)=(1-t)^k h(t)
\end{align}
for some $k\ge 0$ and $h(t)\in \bar K_R[t]$ such that $h(1)\ne 0$.  
Since the leading coefficient of $(1-t)^k$ is a unit in $R$, 
by going through the division of 
$f(t)$ by $(t-1)^k$, it is easy to see that 
$h(t)$ actually lies in $R[t]$.
 
Since $\phi(1)$ is  
an idempotent of $\cA$ and  
$\cA$ is a domain, we have $\phi(1)=0$ or $1$.  
Since $\phi\ne 0$ by assumption,  we have $\phi(1)=1$. 
Applying $0=f(\phi)$ to $1$ we get $f(1)=0$, whence $k\ge 1$. 
Furthermore, since $\phi\ne \I$, we also have $\deg h\ge 1$.

Let $\bar\cA=\bar K_R \otimes_R \cA$. Since $\cA$ is a domain containing $R$, 
and hence torsion-free as an $R$-module, 
the standard map $\cA\simeq R\otimes_R\cA \to K_R\otimes_R \cA$ 
is injective, for by \cite[Prop. 3.3]{AM} 
$K_R\otimes_R \cA$ is isomorphic to the localization $S^{-1}\cA$ 
with $S=R\backslash\{0\}$. Since every field is absolutely flat,  
%(e.g., see \cite{Bou}, pp. 45-46, Ex. 16-17), 
the standard map $K_R\otimes_R\cA \to \bar K_R \otimes_{K_R} (K_R\otimes_R\cA)=
\bar K_R \otimes_R\cA$ is also 
injective. Therefore, we may view $\cA$ as 
an $R$-subalgebra of $\bar\cA$ in the standard way and 
extend $\phi$ $\bar K_R$-linearly to a $\bar K_R$-algebra endomorphism 
of $\bar \cA$, which we will denote by $\bar \phi$. 

Since $\phi$ is algebraic over $R$, $\bar\phi$ is algebraic over 
$\bar K_R$. 
Then $\bar \cA$ can be 
decomposed as a direct sum of the generalized eigen-subspaces of $\bar \phi$ 
(e.g., see  \cite[Proposition $4.2$]{H}).
More precisely, let $r_i$ $(1\le i\le \ell)$ be all the distinct roots of 
$f(t)$ in $\bar K_R$ with multiplicity $m_i$. Set 
$\bar\cA_i=\Ker (r_i\I_{\bar \cA}-\bar\phi)^{m_i}$ for all $1\le i\le \ell$. 
Then we have 
\begin{align}\label{DecomBerA}
\bar \cA=\oplus_{i=1}^\ell \bar \cA_i. 
\end{align} 

Furthermore, the decomposition above is actually an algebra grading of $\bar\cA$, 
i.e., $\bar \cA_i\bar \cA_j\subseteq \bar \cA_{ij}$ for all 
$1\le i, j\le \ell$. In particular, $\bar \cA_1$ is a nonzero 
$R$-subalgebra of $\bar \cA$, and hence also a unital domain over $R$, 
for $1\in \bar\cA_1$. 

Note also that $\bar\cA_1$ is $\phi$-invariant and hence also 
$h(\bar\phi)$-invariant. Furthermore, since $h(1)\ne 0$, 
the restriction of $h(\bar\phi)$ on $\bar \cA_1$ is injective. 
Otherwise, there would exist $0\ne a\in \bar\cA_1$ such that 
$h(\bar\phi)(a)=0$. Since $(\I_{\bar \cA}-\bar\phi)^{m_1}(a)=0$, 
and $h(t)$ and $(1-t)^{m_1}$ are co-prime, we have $a=0$. 
Contradiction.

Now, since $f(\bar\phi)\big |_{\bar\cA_1}=
h(\phi)\big |_{\bar\cA_1}(\I_{\bar\cA_1}-\bar\phi)^k\big |_{\bar\cA_1} =0$, 
we have $(\I_{\bar\cA_1}-\bar\phi)^k\big |_{\bar\cA_1} =0$, i.e., 
$(\I_{\bar\cA_1}-\bar\phi)$ is a nilpotent $R$-$\cE$-derivation of $\bar\cA_1$. 
Then $\I_{\bar\cA_1}-\bar\phi=0$ by Lemma \ref{NoNilHD}. 
Consequently, $\tilde f(\bar\phi)=0$, where $\tilde f(t)=(1-t) h(t)$.
Hence we also have $\tilde f(\phi)=0$. Since $h(t)\in R[t]$ as pointed above, 
we have $\tilde f(t)\in R[t]$. Then by the choice of $f(t)$,  
we have $f(t)=\tilde f(t)$, whence $k=1$, as desired.   
\epfv  
  
\begin{corol}\label{NoLN-ED}
Let $R$ and $\cA$ be as in Lemma \ref{(1-t)Lma}. Then $\cA$ has 
no nonzero locally nilpotent $R$-$\cE$-derivation that 
is algebraic over $R$. 
\end{corol}

\pf Let $\delta\in \Eder_R(\cA)$ be LN and algebraic over $R$.  
Write $\delta=\I-\phi$ for some $\phi\in\cEnd_R(\cA)$. Then 
$\phi=\I-\delta$ is also algebraic over $R$.   
Let $f(t)\in R[t]$ be a minimal polynomial of $\phi$. 
Then for each $a\in \cA$, we have $f(\phi)(a)=0$ and 
$\delta^k(a)=(\I-\phi)^k(a)=0$ for some $k\ge 1$.

Let $K_R$ be the field of fractions of $R$,  $\cB\!:=K_R\otimes_R\cA$, 
and $\bar \phi$ and $\bar\delta$ the $K_R$-linear  
extension maps of $\phi$ and $\delta$, respectively, 
from $\cB$ to $\cB$. As pointed out in the proof of Lemma \ref{(1-t)Lma}, 
we may identify $\cA$ as an $R$-subalgebra of $\cB$. 

With the setting above, we have $\bar\delta=\I_{\cB}-\bar\phi$,  
$f(\bar \phi)(a)=0$ and $\bar \delta^k(a)=(\I_{\bar \cA}-\bar \phi)^k(a)=0$.
By Lemma \ref{(1-t)Lma},  $\gcd (f(t), (1-t)^k)=1-t$ in $K_R[t]$. 
Hence there exist $u(t), v(t)\in K_R[t]$ such that 
$u(t)f(t)+v(t)(1-t)^k=1-t$. Consequently,    
$(\I_\cB-\bar\phi)(a)=0$. Since $a\in\cA$, we further have 
$\delta(a)=(\I-\phi)(a)=(\I_{\cB}-\bar\phi)(a)=0$.  
Therefore, $\delta=0$ and the corollary follows. 
\epfv

From now on we focus on the $\cE$-derivations of integral 
domains of characteristic zero. 

\begin{lemma}\label{phi-Lma}
Assume further that $R$ is an integral domain of characteristic zero, and  
$\cA$ is an integral domain containing $R$.  
Let $\phi\in \cEnd_R(\cA)$ and $g(t)=\sum_{i=r}^d c_i t^i \in R[t]$ with 
$c_r, c_d\ne 0$.   
Then for each $a\in \cA$ such that $a^m\in \Ker g(\phi)$ 
for all $m\ge 1$, the following statements hold:
\begin{enumerate}
  \item[$1)$] $\phi^i(a)=\phi^j(a)$ for some $r\le i<j \le d$; 
  \item [$2)$] if $g(1)\ne 0$, then $\phi^k (a)=0$ for some $r\le k\le d$.  
\end{enumerate} 
\end{lemma}

\pf  If $\phi=0$, the lemma is trivial. 
So we assume $\phi\ne 0$.  If $d=r$, then $g(\phi)=c_r\phi^r$ and 
$\Ker g(\phi)=\Ker \phi^r$. So we have $\phi^r(a)=0$, whence both statements  
$1)$ and $2)$ hold. So we assume $r<d$. 

Set $b_i\!:=\phi^i(a)$ for all $r\le i\le d$. Since $g(\phi) (a^m)=0$ and 
$b_i^m=\phi^i(a)^m=\phi^i(a^m)$ for all $m\ge 1$, we have 
\begin{align}\label{phi-Lma-peq1}
c_rb_r^m +c_{r+1} b_{r+1}^m+ \cdots + c_d b_d^m=0.
\end{align}

Since $\cA$ is an integral domain and not all coefficients $c_i$'s are zero, 
the vandemonde determinant $\prod_{r\le i<j\le d}(b_j-b_i)=0$, whence $b_i-b_j=0$, 
i.e., $\phi^i(a)=\phi^j(a)$, for some $r\le i< j\le d$. So    
statement $1)$ holds.

To show statement $2)$, assume otherwise, i.e., $\phi^i(a)\ne 0$ for all $r\le i\le d$. 
Let $u_k$ $(1\le k\le \ell)$ be all distinct (nonzero) elements of 
$b_i=\phi^i(a)$ $(r \le i \le d)$.
For each $1\le k\le \ell$, let $B_k$ be the subset of $r\le i\le d$ 
such that $\phi^i(a)=u_k$, and set $\tilde c_k\!:=\sum_{i\in B_k} c_i$.
Then $\sum_{k=1}^\ell \tilde c_k=g(1)\ne 0$, whence 
$\tilde c_k$ $(1\le k\le \ell)$ are not all zero.

On the other hand, Eq.\,(\ref{phi-Lma-peq1}) 
above can be re-written as 
\begin{align}\label{phi-Lma-peq2}
\tilde c_1 u_1^m +\tilde c_2 u_2^m+ \cdots + \tilde c_\ell u_\ell^m=0.
\end{align}
Since $\cA$ is an integral domain and $u_k$ $(1\le k\le \ell)$ are 
distinct nonzero elements of $\cA$,
by using the vandemonde determinant we see that $\tilde c_k=0$ 
for all $1\le k\le \ell$. Contradiction.
\epfv

\begin{corol}\label{Corol4.1.9}
Assume that $R$ is an integral domain of characteristic zero, and  
$\cA$ is an integral domain (containing $R$). 
If $\cA$ is finitely generated as an $R$-algebra,  
then for every $\phi\in\cEnd_R(\cA)$ that is algebraic over $R$,  
we have $\phi^i=\phi^j$ for some $1\le i<j$.
\end{corol}

\pf Let $0\ne f(t)\in R[t]$ such that $f(\phi)=0$, and 
$a_k\in\cA $ $(1\le k\le n)$ that generate $\cA$ as 
an $R$-algebra. Hence $\Ker f(\phi)=\cA$ and $a_k\in \rad(\Ker f(\phi))$ for all 
$1\le k\le n$. By lemma \ref{phi-Lma}, for each $1\le k\le n$ 
there exists $1\le i_k<j_k$ such that $\phi^{i_k}(a_k)=\phi^{j_k}(a_k)$.
Applying some powers of $\phi$ to the equation above we may 
assume that $i_k$ $(1\le k\le n)$ are all equal to one another. 
We denote this integer by $i$.   

Set $j=i+\prod_{i=1}^n (j_k-i)$. Then 
it is easy to see that $\phi^i(a_k)=\phi^j(a_k)$ for all $1\le k\le n$. 
Since $\cA$ as an $R$-algebra is generated by 
$a_k$ $(1\le k\le n)$ and $\phi$ is an $R$-algebra 
endomorphism, we have $\phi^i=\phi^j$, as desired.   
\epfv

Now, we are ready to show the main results of this section. 

\begin{propo}\label{HDIntDomain}
Assume that $R$ is an integral domain of characteristic zero,  
$\cA$ an integral domain containing $R$, and $\phi$ an 
$R$-endomorphism of $\cA$ that is algebraic over $R$.  
Set $\Ker_{\ge 1}\phi\!:=\sum_{i\ge 1}\Ker \phi^i$. Then 
we have 
\begin{align}\label{HDIntDomain-eq1}
\rad\big(\im (\I-\phi)\big)=\rad(\Ker_{\ge 1}\phi).
\end{align} 
Consequently, $\im (\I-\phi)$ is a MS of $\cA$. 
\end{propo}

\pf The case $\phi=0$ or $\I$ is trivial, so we assume 
$\phi\ne 0, \I$. By Lemma \ref{ReducLma}, $1)$ we have $ \Ker_{\ge 1}\phi\subseteq 
\im (\I-\phi)$, whence $\rad(\Ker_{\ge 1}\phi)\subseteq \rad(\im (\I-\phi))$. 

Conversely, let $a\in \rad(\im (\I-\phi))$ and $f(t)$ be a minimal polynomial 
of $\phi$. Replacing $a$ by a power of $a$ we assume that 
$a^m\in \im(\I-\phi)$ for all $m\ge 1$. 
Let $K_R$ be the field of fractions of $R$,  
$\cB\!:=K_R\otimes_R\cA$ and $\bar \phi$ the $K_R$-linear 
extension of $\phi$ for $\cB$ to $\cB$. As pointed out  
as in the proof of Lemma \ref{(1-t)Lma}, we may identify $\cA$ as 
an $R$-subalgebra of $\cB$.
%#$ by the flatness of localization prop

By Lemma \ref{(1-t)Lma}, $f(t)=(t-1)h(t)$ for some $h(t)\in R[t]$ 
such that $h(1)\ne 0$ and $\deg h \ge 1$. 
Then there exist $u(t), v(t)\in K_R[t]$ such that 
$(1-t)u(t)+h(t)v(t)=1$. Then by Lemma \ref{GenralLma} with 
$A=\I_{\cB}-\bar \phi$, $B=h(\bar \phi)$, $C=u(\bar \phi)$ 
and $D=v(\bar \phi)$, we have $\im(\I_{\cB}-\bar \phi)=\Ker h(\bar \phi)$, 
whence $h(\phi)(a^m)=h(\bar \phi)(a^m)=0$ for all $m\ge 1$. Applying 
Lemma \ref{phi-Lma}, $2)$ with $g(t)=h(t)$ 
we have $\phi^k(a)=0$ for some $k\ge 0$.  
If $k=0$, then $a=0$, and if $k\ge 1$, $a\in \Ker \phi^k$. 
In either case $a\in \Ker_{\ge 1}\phi$,  
whence Eq.\,(\ref{HDIntDomain-eq1}) follows. 

The statement that $\im (\I-\phi)$ is a MS of $\cA$ follows directly from 
Eq.\,(\ref{HDIntDomain-eq1}), 
Lemma \ref{ReducLma}, $1)$ and Lemma \ref{I-MS-Lma}.
\epfv

\begin{theo}\label{AlgIntDomainThm}
Assume that $R$ is an integral domain of characteristic zero, and  
$\cA$ is an integral domain (containing $R$). Then the LFED conjecture  
(resp., the LNED conjecture) holds for all   
(resp., locally nilpotent) $R$-derivations and 
$R$-$\cE$-derivations of $\cA$ that are algebraic over $R$. 
\end{theo} 

\pf By Proposition \ref{NoAlgDer}, $\cA$ has no nonzero  
%LN (locally nilpotent) 
$R$-derivation that is algebraic over $R$. 
Hence the $R$-derivation case of the corollary holds.

By Corollary \ref{NoLN-ED}, $\cA$ has no nonzero  
locally nilpotent $R$-$\cE$-derivation that is algebraic over $R$. 
Hence the $R$-$\cE$-derivation case of the LNED conjecture 
in the corollary holds. The $R$-$\cE$-derivation case of 
the LFED conjecture in the corollary follows directly 
from Proposition \ref{HDIntDomain}.
\epfv

\end{document}